# The equivariant moving frame method of third order differential equations

*M. Nadjafikhah*[†] *and R. Bakhshandeh Chamazkoti*[‡]

[†]*School of Mathematics, Iran University of Science and Technology, Narmak, Tehran 1684613114, Iran.*

[‡]*Department of Mathematics, Faculty of Basic Science, Babol University of Technology, Babol, Iran*

E-mail: [†]*m_nadjafikhah@iust.ac.ir*, [‡]*r_bakhshandeh@nit.ac.ir*

**Abstract.** This paper is devoted to apply the equivariant moving frame method to study the local equivalence problem of third order ordinarily differential equation under the pseudo-group of fiber preserving transformations.

## 1 Introduction

The general equivalence problem is to recognize when two geometrical objects are mapped on each other by a certain class of diffeomorphisms. E. Cartan developed the general equivalence problem and provided a systematic procedure for determining the necessary and sufficient condition [1, 2, 3]. Now we assume $y = y(x)$ be a real function and give two third order ordinary differential equation $y''' = F(x, y, y', y'')$ and $\bar{y}''' = \bar{F}(\bar{x}, \bar{y}, \bar{y}', \bar{y}'')$. The equivalence problem is to establish whether or not there exists a local transformation of variables of a suitable type that transforms first into second equation. This problem involves a number of related problems such as defining a class of transformations, finding invariants of these transformations, obtaining the equivalence criteria and constructing the transformations. The main contribution to the issue of point and contact geometry of third order ODEs was made by respectively E. Cartan and S. -S. Chern in their classical papers [1, 2, 3] and [4]. One may consider equivalence with respect to several types of transformations, for example Michal Godlinski was extensively studied on contact, point and fibre-preserving transformations cases for third order ODEs, in his Ph.D. thesis in [5].

But in this attempt we shall apply the equivariant moving frame method to solve the local equivalence problem of third order differential equations under the pseudo-group of fiber preserving transformations

$$X = \xi(x), \qquad U = \varphi(x, u). \tag{1.1}$$

The theory of moving frames was primarily developed and formulated by Élie Cartan, [6, 7], who modeled it into an algorithmic tool for determining a generating set of the differential invariant algebra of Lie pseudo-groups. The moving frame provides an effective means for determining complete systems of differential invariants and invariant differential forms, classifying their syzygies and recurrence relations, and solving equivalence and symmetry problems arising in a broad range of applications, [8].

In [9], Fels and Olver suggested a new theoretical foundation to the method of moving frames. Suppose a Lie group $G$ acting on the $n$–th order jet space $J^n$ of $m$-dimensional submanifolds of $M$ is given. A moving frame is a $G$–equivariant section of the trivial bundle $J^n \times G \longrightarrow J^n$. This new concept of moving frames is now referred as the *equivariant moving frame method*. In the finite-dimensional theory, [9], a moving frame is defined as an equivariant map $J^n \longrightarrow G$ from an open subset of the submanifold jet bundle to the Lie group. For Lie pseudo-groups, we define a moving frame to be an equivariant section of a suitable bundle $\mathcal{H}^{(n)} \longrightarrow J^n$ constructed from the jets of pseudo-group transformations. For finite-dimensional Lie group actions, the existence of a moving frame requires that the action be free, i.e., $G_S = \{g \in G \ : \ g \cdot S = S\}$ which denote the *isotropy subgroup* of a subset $S \subset M$ be trivial. Clearly, an infinite-dimensional pseudo-group action never has trivial isotropy, and so we must modify the definition of freeness to require that all elements of the isotropy sub-pseudo-group of a point in $J^n$ have the same $n$–th order jet as the identity diffeomorphism. Indeed, Freeness of a prolonged pseudo-group action does not reduce to the usual freeness condition when the pseudo-group is a finite-dimensional Lie group.

Assuming freeness, one can construct the moving frame using Cartan normalization procedure associated with a choice of local cross-section to the group orbits in $J^n$, [9, 10]. The moving frame induces an invariantization process that canonically maps general differential functions and differential forms on $J^\infty$ to their invariant counterparts.

New and significant applications of the equivariant moving frame method for infinite-dimensional Lie pseudogroups actions appears in [11, 12, 13] by Olver and Pohjanpelto. There are some applications of this theory for infinite-dimensional Lie pseudo-groups in [14, 15].

After reviewing this theory, the method is illustrated for third order differential equations. The theory of infinite-dimensional Lie pseudo-groups relies on the Cartan–Kähler Theorem, [16, 17], which requires analyticity. Thus, all our constructions and results hold in the analytic category.

## 2   Structure Equations of Lie Pseudo-Groups

Throughout this paper, $M$ will be an $m$-dimensional manifold and $\mathcal{D} = \mathcal{D}(M)$ the pseudo-group of all local analytic diffeomorphisms of $M$. A pseudo-group $\mathcal{G} \subset \mathcal{D}$ is called *regular* of order $n^\star \geq 1$ if, for all finite $n \geq n^\star$, the $n$–jet of the transformations form an embedded subbundle $\mathcal{G}^{(n)} \subset \mathcal{D}^{(n)}$, and the projection $\pi_n^{n+1} : \mathcal{G}^{(n+1)} \to \mathcal{G}^{(n)}$ is a surjective submersion. An analytic pseudo- group $\mathcal{G} \subset \mathcal{D}$ is called a *Lie pseudo-group* if $\mathcal{G}$ is regular of order $n^\star \geq 1$ and, moreover, every local diffeomorphism $\varphi$ of $\mathcal{D}$ satisfying $j_{n^\star}\varphi \subset \mathcal{G}^{(n^\star)}$ belongs to the pseudo-group, i.e. $\varphi \in \mathcal{G}$.

In local coordinates, for $n \geq n^\star$, the pseudo-group jet subbundle $\mathcal{G}^{(n)} \subset \mathcal{D}^{(n)}$ is characterized by an involutive system of $n^{\text{th}}$ order partial differential equations

$$F^{(n)}(z, Z^{(n)}) = 0, \tag{2.1}$$

called the *determining system* of $\mathcal{G}$, and a local diffeomorphism $\varphi$ is in $\mathcal{G}$ if and only if its $n$–jet is a solution of (2.1).

In local coordinates a vector field in $\mathcal{X}(M)$, the space of locally defined vector fields on $M$, is denote by

$$\mathbf{v} = \sum_{a=1}^m \zeta^a(z) \frac{\partial}{\partial z^a}, \tag{2.2}$$

For $0 \leq n \leq \infty$, let $J^n TM$ denote the bundle of $n$th order jets of sections of $TM$. Local coordinates on $J^n TM$ are given by

$$(z, \zeta^{(n)}) = (z^a, \zeta_B^a), \qquad a = 1, \ldots, m, \qquad 0 \leq \#B \leq n, \tag{2.3}$$

where $\zeta_B^a$ represents the partial derivative $\partial^{\#B}\zeta^a/\partial z^B$. Let

$$L^{(n)}(z,\zeta^{(n)}) = \sum_{b=1}^m \sum_{\#A \leq n} h_b^A(z)\zeta_A^b = 0, \tag{2.4}$$

be the linear system of partial differential equations obtained by linearizing the determining system (2.1) at the identity $n$-jet $\mathbf{1}^{(n)} \in \mathcal{G}^{(n)}$. A vector field (2.2) is in the Lie algebra $\mathfrak{g}$ of infinitesimal generators of the Lie pseudo-group $\mathcal{G}$ if and only if its $n$-jet is a solution of (2.4).

The structure equations of an infinite-dimensional Lie pseudo-group $\mathcal{G}$ are obtained by computing the structure equations of a $\mathcal{G}$-invariant coframe. For the diffeomorphism pseudo-group $\mathcal{D}$ an invariant coframe is given by the horizontal forms

$$\sigma^a = d_M Z^a = \sum_{b=1}^m Z_b^a dz^b, \qquad a = 1,\ldots,m, \tag{2.5}$$

while

$$\mu^a = d_G Z^a = dZ^a - \sum_{b=1}^m z_b^a dz^b, \qquad b = 1,\ldots,m, \tag{2.6}$$

are the zero–th order invariant contact forms and the Maurer–Cartan forms

$$\mu_A^b, \quad b = 1,\ldots,m, \ \#A \geq 0, \tag{2.7}$$

on $\mathcal{D}^\infty$. Writing the horizontal component of the exterior differential of a differential function $F : \mathcal{D}^{(\infty)} \to \mathbb{R}$ in terms of the invariant horizontal coframe (2.5)

$$d_M F = \sum_{a=1}^m (\mathbb{D}_{Z^a} F)\sigma^a, \tag{2.8}$$

serves to define the dual invariant total differential equations

$$\mathbb{D}_{Z^a} = \sum_{b=1}^m w_a^b \mathbb{D}_{z^b}, \qquad a = 1,\ldots,m, \tag{2.9}$$

where $(w_a^b(z,Z^{(1)})) = (\partial Z^b/\partial z^a)^{-1}$ denotes the inverse Jacobian matrix.

**Theorem 2.1.** *The structure equations of the diffeomorphism pseudo-group $\mathcal{D}(M)$ are*

$$d\sigma = \sum_{b=1}^m \mu_b^a \wedge \sigma^b, \quad d\mu_C^a = \sum_{b=1}^m \left[\sigma^b \wedge \mu_{C,b}^a + \sum_{\substack{C=(A,B) \\ \#B \geq 1}} \binom{C}{A} \mu_{A,b}^a \wedge \mu_B^b \right], \tag{2.10}$$

*where the last sum ranges over all multi-indices $A, B$ such that $\#B \geq 1$, and their concatenation equals the multi-index $C$ where*

$$\binom{C}{A} = \frac{C!}{A!B!}, \quad \text{when} \quad C = (A,B).$$

**Proposition 2.2.** *Let $\mathcal{G} \subset \mathcal{D}$ be a Lie pseudo-group. Once restricted to $\mathcal{G}^\infty$ the Maurer–Cartan forms (2.7) satisfy the linear relations*

$$L^{(\infty)}(Z,\mu^{(\infty)}) = 0, \tag{2.11}$$

*obtained by making the formal replacements $z^a \mapsto Z^a$, $\zeta_B^a \mapsto \mu_B^a$ in the infinitesimal determining equations (2.4).*

**Theorem 2.3.** *The structure equations of a Lie pseudo-group $\mathcal{G} \subset \mathcal{D}$ are obtained by restricting the structure equations (2.10) of the diffeomorphism pseudo-group to the kernel of (2.11).*

We compute the structure equations of the Lie pseudo-group action

$$\begin{aligned}
X &= \xi(x), \qquad U = \varphi(x,u), \qquad P = (p\,\varphi_u + \varphi_x)/\xi_x, \\
Q &= \left[\xi_x \varphi_{uu} p^2 - (\xi_{xx}\varphi_u - 2\xi_x\varphi_{xu})p + \varphi_u \xi_x q - \xi_{xx}\varphi_x + \xi_x\varphi_{xx}\right]/\xi_x^3, \\
R &= -\big[-\xi_x^2 \varphi_{uuu} p^3 + 3\xi_x(\xi_{xx}\varphi_{uu} - \xi_x\varphi_{xuu})p^2 + (6\xi_{xx}\xi_x\varphi_{xu} + \xi_{xxx}\xi_x\varphi_u - 3\xi_{xx}^2\varphi_u \\
&\quad -3\xi_x^2 \varphi_{xxu})p + 3\xi_x(\varphi_u \xi_{xx} - \varphi_{xu}\xi_x)q - 3\xi_x^2 \varphi_{uu} pq - \varphi_u \xi_x^2 r + 3\xi_{xx}\xi_x\varphi_{xx} + \xi_{xxx}\xi_x\varphi_x \\
&\quad -\xi_x^2 \varphi_{xxx} - 3\xi_{xx}^2\varphi_x\big]/\xi_x^5,
\end{aligned} \qquad (2.12)$$

which is obtained from (1.1) by considering its third order prolongation. The determining equations for the Lie pseudo-group action (2.12) are

$$X_u = X_p = X_q = X_r = 0, \quad U_p = U_q = U_r = 0, \quad P = (p\,U_u + U_x)/X_x. \qquad (2.13)$$

To obtain the infinitesimal determining equations of the pseudo-group we linearize (2.12) at the identity transformation. If an infinitesimal generator of the pseudo-group action is denoted by

$$\mathbf{v} = \alpha(x,u,p,q,r)\frac{\partial}{\partial x} + \beta(x,u,p,q,r)\frac{\partial}{\partial u} + \gamma(x,u,p,q,r)\frac{\partial}{\partial p} + \tau(x,u,p,q,r)\frac{\partial}{\partial q}$$
$$+ \varsigma(x,u,p,q,r)\frac{\partial}{\partial r}, \qquad (2.14)$$

then the infinitesimal determining system is given by

$$\begin{aligned}
&\alpha_u = \alpha_p = \alpha_q = \alpha_r = 0, \quad \beta_p = \beta_q = \beta_r = 0, \\
&\gamma = (\beta_u - \alpha_x)p + \beta_x, \\
&\tau = \beta_{uu}p^2 + (2\beta_{xu} - \alpha_{xx})p + (\beta_u - 2\alpha_x)q + \beta_{xx}, \\
&\varsigma = \beta_{uuu}p^3 + 3\beta_{xuu}p^2 + (3\beta_{xxu} - \alpha_{xxx})p + 3(\beta_{xu} - \alpha_{xx})q + (\beta_u - 3\alpha_x)r + 3\beta_{uu}pq + \beta_{xxx}.
\end{aligned} \qquad (2.15)$$

Using the propositionosition 2.1 we find the linear relations among the Maurer-Cartan forms $\mu_A^x$, $\mu_A^u$, $\mu_A^p$, $\mu_A^q$, $\mu_A^r$ and their prolongations by the following substituting

$$\alpha_A \mapsto \mu_A^x, \quad \beta_A \mapsto \mu_A^u, \quad \eta_A \mapsto \mu_A^p, \quad \tau_A \mapsto \mu_A^q, \quad \varsigma_A \mapsto \mu_A^r, \qquad (2.16)$$

in the infinitesimal determining equations (2.15) and then the lift of it gives the linear relations

$$\begin{aligned}
&\mu_U^x = \mu_P^x = \mu_Q^x = \mu_R^x = 0 \qquad \mu_P^u = \mu_Q^u = \mu_R^u = 0, \\
&\mu^p = P(\mu_U^u - \mu_X^x) + \mu_X^u, \\
&\mu^q = \mu_{UU}^u P^2 + (2\mu_{XU}^u - \mu_{XX}^x)P + (\mu_U^u - 2\mu_X^x)Q + \mu_{XX}^u, \\
&\mu^r = \mu_{UUU}^u P^3 + 3\mu_{XUU}^u P^2 + (3\mu_{XXU}^u - \mu_{XXX}^x)P + 3(\mu_{XU}^u - \mu_{XX}^x)Q \\
&\quad + (\mu_U^u - 3\mu_X^x)R + 3\mu_{UU}^u PQ + \mu_{XXX}^u.
\end{aligned} \qquad (2.17)$$

It follows that

$$\mu_i = \mu_{X^i}^x, \quad \mu_{i,j} = \mu_{U^i X^j}^u, \quad i,j \geq 0, \qquad (2.18)$$

is a basis of Maurer–Cartan forms for the pseudo-group action (2.12). To simplify the calculation, we introduce the functions

$$\delta(x,u,p) = p\varphi_u + \varphi_x, \qquad \varepsilon(x) = \xi_{xx}/\xi_x,$$
$$\psi(x,u,p,q) = \xi_x\varphi_{uu}p^2 - (\xi_{xx}\varphi_u - 2\xi_x\varphi_{xu})p + \varphi_u\xi_x q - \xi_{xx}\varphi_x + \xi_x\varphi_{xx}, \qquad (2.19)$$
$$\chi(x,u,p,q,r) = -\xi_x^2\varphi_{uuu}p^3 + 3\xi_x(\xi_{xx}\varphi_{uu} - \xi_x\varphi_{uu})p^2 + (6\xi_{xx}\xi_x\varphi_{xu} + \xi_{xxx}\xi_x\varphi_u - 3\xi_{xx}^2\varphi_u$$
$$-3\xi_x^2\varphi_{xxu})p + 3\xi_x(\varphi_u\xi_{xx} - \varphi_{xu}\xi_x)q - 3\xi_x^2\varphi_{uu}pq - \varphi_u\xi_x^2 r + 3\xi_{xx}\xi_x\varphi_{xx} + \xi_{xxx}\xi_x\varphi_x$$
$$-\xi_x^2\varphi_{xxx} - 3\xi_{xx}^2\varphi_x.$$

The horizontal forms of the pseudo-group action (2.12) are given by

$$\sigma^x = \xi_x dx, \qquad \sigma^u = \varphi_x dx + \varphi_u du,$$
$$\sigma^p = \frac{1}{\xi_x}[(\delta_x - \delta\varepsilon)dx + \delta_u du + \varphi_u dp],$$
$$\sigma^q = \frac{1}{\xi_x^3}[(\psi_x - 3\psi\varepsilon)dx + \psi_u du + \psi_p dp + \xi_x\varphi_u dq], \qquad (2.20)$$
$$\sigma^r = -\frac{1}{\xi_x^5}[(\chi_x - 5\chi\varepsilon)dx + \chi_u du + \chi_p dp + \chi_q dq - \varphi_u\xi_x^2 dr]$$

Applying Theorem 2.3, the structure equations of the invariant horizontal coframe are

$$d\sigma^x = -d\mu^x = \mu_X^x \wedge \sigma^x,$$
$$d\sigma^u = -d\mu^u = \mu_X^u \wedge \sigma^x + \mu_U^u \wedge \sigma^u,$$
$$d\sigma^p = -d\mu^p = \mu_X^p \wedge \sigma^x + \mu_U^p \wedge \sigma^u + \mu_P^p \wedge \sigma^p, \qquad (2.21)$$
$$d\sigma^q = -d\mu^q = \mu_X^q \wedge \sigma^x + \mu_U^q \wedge \sigma^u + \mu_P^q \wedge \sigma^p + \mu_Q^q \wedge \sigma^q,$$
$$d\sigma^r = -d\mu^r = \mu_X^r \wedge \sigma^x + \mu_U^r \wedge \sigma^u + \mu_P^r \wedge \sigma^p + \mu_Q^r \wedge \sigma^q + \mu_R^r \wedge \sigma^r.$$

where

$$\mu_X^p = P(\mu_{UX}^u - \mu_{XX}^x) + \mu_{XX}^u, \qquad \mu_U^p = P\mu_{UU}^u + \mu_{XU}^u, \qquad \mu_P^p = \mu_U^u - \mu_X^x,$$
$$\mu_X^q = \mu_{UUX}^u P^2 + (2\mu_{XXU}^u - \mu_{XXX}^x)P + (\mu_{UX}^u - 2\mu_{XX}^x)Q + \mu_{XXX}^u,$$
$$\mu_U^q = \mu_{UUU}^u P^2 + 2\mu_{XUU}^u P + \mu_{UU}^u Q + \mu_{XXU}^u,$$
$$\mu_P^q = 2\mu_{UU}^u P + 2\mu_{XU}^u - \mu_{XX}^x, \qquad \mu_Q^q = \mu_U^u - 2\mu_{XX}^x,$$
$$\mu_X^r = \mu_{UUUX}^u P^3 + 3\mu_{XXUU}^u P^2 + (3\mu_{XXXU}^u - \mu_{XXXX}^x)P + 3(\mu_{XXU}^u - \mu_{XXX}^x)Q \qquad (2.22)$$
$$+ (\mu_{UX}^u - 3\mu_{XX}^x)R + 3\mu_{UUX}^u PQ + \mu_{XXXX}^u,$$
$$\mu_U^r = \mu_{UUUU}^u P^3 + 3\mu_{XUUU}^u P^2 + 3\mu_{XXUU}^u P + 3\mu_{XUU}^u Q + \mu_{UU}^u R + 3\mu_{UUU}^u PQ + \mu_{XXXU}^u,$$
$$\mu_P^r = 3\mu_{UUU}^u P^2 + 6\mu_{XUU}^u P + 3\mu_{XXU}^u - \mu_{XXX}^x + 3\mu_{UU}^u Q,$$
$$\mu_Q^r = 3(\mu_{UU}^u P + \mu_{XU}^u - \mu_{XX}^x), \qquad \mu_R^r = \mu_U^u - 3\mu_X^x.$$

Also, the structure equations for the basis of Maurer-Cartan forms (2.18) are

$$d\mu_i = \sigma^x \wedge \mu_{i+1} + \sum_{0 \le a < i}\binom{i}{a}\mu_{a+1} \wedge \mu_{i-a},$$
$$d\mu_{i,j} = \sigma^x \wedge \mu_{i,j+1} + \sigma^u \wedge \mu_{i+1,j} + \sum_{0 \le a < j}\binom{j}{a}\mu_{i,a+1} \wedge \mu_{j-a}. \qquad (2.23)$$

## 3 Equivariant Moving Frames

For infinite-dimensional Lie pseudo-group actions, the equivariant moving frame construction was first laid out in [11].

Let $\mathcal{G}$ be a Lie pseudo-group action on an $m$-dimensional manifold $M$ and $1 \leq p < m$. For each integer $0 \leq n \leq \infty$, let $J^n = J^n(M,p)$ denote the the extended jet bundle of equivalence classes of p-dimensional submanifolds of $M$ under $n^{\text{th}}$ order contact at a single point, [17]. For $k \geq n$, we use $\pi_n^k : J^k \to J^n$ to denote the canonical projection. Locally, the coordinates on $M$ can be written as $z = (x,u)$ where $x = (x^1, \ldots, x^p)$ are considered to be the independent variables parameterizing a submanifold $S \subset M$ and $u = (u^1, \ldots, u^q)$, $q = m - p$, the dependent variables. The induced coordinates on $J^n$ are denoted by $z^{(n)} = (x, u^{(n)})$, where $z^{(n)}$ denotes the derivatives $u_J^\alpha = \partial^{\#J} u^\alpha / \partial x^J$ of the $u$'s with respect to the $x$'s of order $0 \leq \#J \leq n$.

Let $\mathcal{E}^{(n)} \to J^n$ be the lifted bundle obtained by taking the pull-back bundle of $\mathcal{E}^{(n)} \to M$ via the projection $\pi_0^n : J^n \to M$. Local coordinates on $\mathcal{E}^{(n)}$ are given by $(z^{(n)}, g^{(n)})$, where the base coordinates $z^{(n)} = (x, u^{(n)}) \in J^n$ are the submanifold jet coordinates and the fiber coordinates $g^{(n)}$ parameterize the pseudo-group jets. The bundle $\mathcal{E}^{(n)}$ carries the structure of a groupoid, with source map $\sigma(z^{(n)}, g^{(n)}) = z^{(n)}$ and target map $\tau(z^{(n)}, g^{(n)}) = Z^{(n)} = g^{(n)} \cdot z^{(n)}$ given by the prolonged action. The local coordinate expressions for the prolonged action $Z^{(n)}$ are obtained by implementing the chain rule. Suppose that the

$$d_H X^i = \sum_{j=1}^p (D_{x^j} X^i) dx^j, \qquad i = 1, \ldots, p, \tag{3.1}$$

be the lifted horizontal coframe on $\mathcal{E}^{(\infty)}$, where

$$D_{x^j} = \frac{\partial}{\partial x^j} + \sum_{\alpha=1}^q \sum_{\#J \geq 0} u_{J,j}^\alpha \frac{\partial}{\partial u_J^\alpha}, \qquad j = 1, \ldots, p, \tag{3.2}$$

are the *total derivative operators* on the submanifold jet bundle $J^\infty$. The *lifted total differential equations* are defined by the formula

$$d_H F(z^{(n)}) = \sum_{i=1}^p (D_{x^i}) dx^i = \sum_{i=1}^p (D_{X^i} F) d_H X^i. \tag{3.3}$$

where

$$D_{X^i} = \sum_{i=1}^p W_i^j D_{x^i}, \qquad (W_i^j) = (D_{x^i} X^j)^{-1}, \tag{3.4}$$

Let $\mathcal{G}$ be a regular Lie pseudo-group acting on $M$. An $n$–th order *moving frame* is a $\mathcal{G}$-equivariant local section $\rho^{(n)} : J^n \to \mathcal{E}^{(n)}$. A moving frame exists in a neighborhood of a jet $z^{(n)}$ if the action is *free* and *regular*. Let $\mathcal{V}^n$ be the set of regular $n$–th order submanifold jets. A moving frame is obtained by choosing a cross-section $\mathcal{K}^n$ to the pseudo-group orbits and normalizing the pseudo-group parameters Maurer-Cartan forms.

To obtain the prolonged action we apply the lifted total differential equations

$$D_R = \frac{\xi_x^3}{\varphi_u} D_r, \qquad D_Q = \frac{1}{\varphi_u}\Big[\xi_x^2 D_q + \frac{\chi_q}{\xi_x^3} D_R\Big],$$

$$D_P = \frac{1}{\varphi_u}\Big[\xi_x D_p - \frac{\psi_p}{\xi_x^2} D_Q + \frac{\chi_p}{\xi_x^4} D_R\Big],$$

$$D_U = \frac{1}{\varphi_u}\Big[D_u - \frac{\delta_u}{\xi_x} D_P - \frac{\psi_u}{\xi_x^3} D_Q + \frac{\chi_u}{\xi_x^5} D_R\Big],$$

$$D_X = \frac{1}{\xi_x}\Big[D_x - \varphi_x D_U - \big(\frac{\delta_x - \delta\varepsilon}{\xi_x}\big) D_P - \big(\frac{\psi_x - 3\psi\varepsilon}{\xi_x^3}\big) D_Q + \big(\frac{\chi_x - 5\chi\varepsilon}{\xi_x^5}\big) D_R\Big], \qquad (3.5)$$

to the lifted invariant (2.12). For instance, by applying the lifted total differential equations (3.5) to the lifted invariant $R$ yields the first order lifted differential invariants as follows

$$R_Q = -\frac{1}{\varphi_u \xi_x^5}\Big[3\xi_x^3(\varphi_u \xi_{xx} - \varphi_{xu}\xi_x - \xi_x^2 \varphi_{uu} p) - \varphi_u \xi_x^2 r_q + \frac{\chi_q}{\xi_x^3}\Big],$$

$$R_P = -\frac{1}{\varphi_u \xi_x^4}\Big[-\xi_x^2(3\varphi_{uuu}p^2 + 3\varphi_{xxu} + 3\varphi_{uu}q + \varphi_u r_p) + 6\xi_x(\xi_{xx}\varphi_{uu} - \xi_x \varphi_{uu})p$$

$$-3\xi_{xx}^2\varphi_u + 6\xi_{xx}\xi_x\varphi_{xu} + \xi_{xxx}\xi_x\varphi_u - \frac{\psi_p}{\xi_x^3} R_Q + \frac{\chi_p}{\xi_x^5}\Big], \qquad (3.6)$$

$$R_U = -\frac{1}{\varphi_u \xi_x^5}\Big[\chi_u - \frac{\delta_u}{\xi_x} R_P - \frac{\psi_u}{\xi_x^3} R_Q + \frac{\chi_u}{\xi_x^5}\Big],$$

$$R_X = -\frac{1}{\xi_x^6}\Big[\chi_x - \varphi_x R_U - \big(\frac{\delta_x - \delta\varepsilon}{\xi_x}\big) R_P - \big(\frac{\psi_x - 3\psi\varepsilon}{\xi_x^3}\big) R_Q + \frac{\chi_x - 5\chi\varepsilon}{\xi_x^5}\Big],$$

Differentiating (3.6) with respect to (3.5) yields the higher order lifted differential invariants.

Once a moving frame is obtained it is possible to systematically invariantize differential functions, differential forms and differential equations. The space of differential forms on $\mathcal{E}^{(\infty)}$ splits into

$$\mathbf{\Omega}^* = \bigoplus_{k,l} \mathbf{\Omega}^{k,l} = \bigoplus_{i,j,l} \mathbf{\Omega}^{i,j,l}, \qquad (3.7)$$

where $l$ indicates the number of Maurer–Cartan forms, and $k = i + j$ the number of *jet forms*, with $i$ indicating the number of horizontal forms $dx^i, 1 \le i \le p$, and $j$ the number of basic contact forms

$$\theta_J^\alpha = du_J^\alpha - \sum_{i=1}^p u_{J,i}^\alpha dx^i, \qquad \alpha = 1, \ldots, q, \qquad \#J \ge 0, \qquad (3.8)$$

on the extended submanifold jet bundle $J^\infty$. Let

$$\mathbf{\Omega}_J^* = \bigoplus_k \mathbf{\Omega}^{k,0} = \bigoplus_{i,j} \mathbf{\Omega}^{i,j,0}, \qquad (3.9)$$

denote the subspace of jet forms consisting of those differential forms containing no Maurer–Cartan forms. Let $\pi_J : \mathbf{\Omega}^* \to \mathbf{\Omega}_J^*$ be the natural projection that takes a differential form $\Omega$ on $\mathcal{E}^\infty$ to its jet component $\pi_J(\Omega)$ obtained by annihilating all Maurer–Cartan forms in $\Omega$.

The lift of a differential form $\omega$ on $J^\infty$ is the jet form

$$\lambda(\omega) = \pi_J[\tau^*(\omega)]. \qquad (3.10)$$

The lift of a vector jet coordinate $\zeta_A^b$ is defined to be the Maurer–Cartan form $\mu_A^b$

$$\lambda(\zeta_A^b) = \mu_A^b \qquad b = 1,\ldots,m, \qquad \#A \geq 0.$$

**Theorem 3.1.** *Let $\omega$ be a differential form on $J^\infty$. Then*

$$d[\lambda(\omega)] = \lambda[d\omega + \mathbf{v}^{(\infty)}(\omega)], \tag{3.11}$$

*where $\mathbf{v}^{(\infty)}$ is the vector field*

$$\mathbf{v}^{(\infty)} = \sum_{i=1}^p \xi^i(x,u)\frac{\partial}{\partial x^i} + \sum_{\alpha=1}^q \sum_{k=\#J\geq 0} \phi^{\alpha;J}(x,u^{(k)})\frac{\partial}{\partial u_J^\alpha} \in \mathfrak{g}^{(\infty)}, \tag{3.12}$$

*with $\phi^{\alpha;J}$ defined recursively by the prolongation formula, [17],*

$$\phi^{\alpha;J,j} = D_{x^j}\phi^{\alpha;J} - \sum_{i=1}^p (D_{x^j}\xi^i)u_{J,i}^\alpha. \tag{3.13}$$

In particular, the identity (3.11) applies to the lifted differential invariants $X^i, \widehat{U}_J^\alpha$ leading to

$$dX^i = \omega^i + \mu^i, \quad 1 \leq i \leq p, \qquad d\widehat{U}_J^\alpha = \sum_{j=1}^p \widehat{U}_{J,j}^\alpha \omega^j + \widehat{\phi}^{\alpha;J}, \quad 1 \leq \alpha \leq q, \quad \#J \geq 0, \tag{3.14}$$

where $\omega^i = \lambda(dx^i)$ for $i = 1,\ldots,p$ are the lifts of the jet forms, and $\widehat{\phi}^{\alpha;J} = \lambda(\phi^{\alpha;J})$ are *correction terms* obtained by lifting the prolonged vector field coefficients (3.13).

In our problem, we have the equalities

$$\omega^x = \lambda(dx) = \sigma^x, \qquad \omega^u = \lambda(du) = \sigma^u, \qquad \omega^p = \lambda(dp) = \sigma^p,$$
$$\omega^q = \lambda(dq) = \sigma^q, \qquad \omega^r = \lambda(dr) = \sigma^r, \tag{3.15}$$

where the coordinate expressions for the horizontal forms are given in (2.11). To obtain the correction terms in the recurrence relations (3.14) one needs to compute the prolongation of the infinitesimal generator

$$\mathbf{v} = \alpha(x)\frac{\partial}{\partial x} + \beta(x,u)\frac{\partial}{\partial u} + [(\beta_u - \alpha_x)p + \beta_x]\frac{\partial}{\partial p} + [\beta_{uu}p^2 + (2\beta_{xu} - \alpha_{xx})p + (\beta_u - 2\alpha_x)q$$
$$+ \beta_{xx}]\frac{\partial}{\partial q} + [\beta_{uuu}p^3 + 3\beta_{xuu}p^2 + (3\beta_{xxu} - \alpha_{xxx})p + 3(\beta_{xu} - \alpha_{xx})q + (\beta_u - 3\alpha_x)r$$
$$+ 3\beta_{uu}pq + \beta_{xxx}]\frac{\partial}{\partial r}, \tag{3.16}$$

using the formula (3.13). Now we can use the recurrence relations (3.14) to investigate the local equivalence problem for third order differential equations under the group of transformations (2.2). For this aim, firstly we construct a moving frame. This is done by determining a cross-section to the equivalence pseudo-group action. Since there is a correspondence between the normalization of the pseudo-group parameters and the normalization of the Maurer–Cartan forms (2.9), the recurrence relations (3.14) can be used to find a cross-section. First, We restrict our attention to the zero-order lifted differential invariants and obtain, using (3.14), the recurrence relations

$$d_{\mathcal{G}}X \equiv \mu^x, \qquad d_{\mathcal{G}}U \equiv \mu^u, \qquad d_{\mathcal{G}}P \equiv P(\mu_U^u - \mu_X^x) + \mu_X^u,$$
$$d_{\mathcal{G}}Q \equiv \mu_{UU}^u P^2 + (2\mu_{XU}^u - \mu_{XX}^x)P + (\mu_U^u - 2\mu_X^x)Q + \mu_{XX}^u, \tag{3.17}$$
$$d_{\mathcal{G}}R \equiv \mu_{UUU}^u P^3 + 3\mu_{XUU}^u P^2 + (3\mu_{XXU}^u - \mu_{XXX}^x)P + 3(\mu_{XU}^u - \mu_{XX}^x)Q$$
$$+ (\mu_U^u - 3\mu_X^x)R + 3\mu_{UU}^u PQ + \mu_{XXX}^u,$$

where $d_\mathcal{G}$ equals to $\pi_\mathcal{G} \circ d$. Since the group differential of $X, U, P, Q, R$ depend on the linearly independent Maurer–Cartan forms $\mu^x, \mu^u, \mu^u_X, \mu^u_{XX}, \mu^u_{XXX}$ respectively, it is possible to translate $X, U, P, Q$ and $R$ to zero. Setting

$$X = U = P = Q = R = 0, \tag{3.18}$$

leads to the normalizations

$$\omega^x \equiv -\mu^x, \quad \omega^u \equiv -\mu^u, \quad \omega^p \equiv -\mu^u_X, \quad \omega^q \equiv -\mu^u_{XX}. \tag{3.19}$$

The group differential of the first order lifted invariants reduce to

$$\begin{aligned}
d_\mathcal{G} R_Q &\equiv 3(\mu^u_{XU} - \mu^x_{XX}) + \mu^x_X R_Q + 3\mu^u_{UU} P, \\
d_\mathcal{G} R_P &\equiv 3\mu^u_{UUU} P^2 + 6\mu^u_{XUU} P + 3\mu^u_{XXU} - \mu^x_{XXX} + 3\mu^u_{UU} Q - [2\mu^u_{UU} P + 2\mu^u_{XU} \\
&\quad -\mu^x_{XX}] R_Q - 2\mu^x_X R_P, \\
d_\mathcal{G} R_U &\equiv \mu^u_{UUUU} P^3 + 3\mu^u_{XUUU} P^2 + 3\mu^u_{XXUU} P + 3\mu^u_{XUU} Q + \mu^u_{UU} R - 3\mu^x_X R_U \\
&\quad + 3\mu^u_{UUU} PQ - (P\mu^u_{UU} + \mu^u_{XU}) R_P - [\mu^u_{UUU} P^2 + 2\mu^u_{XUU} P + \mu^u_{UU} Q] R_Q, \\
d_\mathcal{G} R_X &\equiv \mu^u_{XUUU} P^3 + 3\mu^u_{XXUU} P^2 + (3\mu^u_{XXXU} - \mu^x_{XXXX}) P + 3(\mu^u_{XXU} - \mu^x_{XXX}) Q \\
&\quad + (\mu^u_{XU} - 3\mu^x_{XX}) R + 3\mu^u_{XUU} PQ + \mu^u_{XXXX} + (\mu^u_U - 4\mu^x_X) R_X - \mu^u_X R_U \\
&\quad - [P(\mu^u_{XU} - \mu^x_{XX}) + \mu^x_{XX}] R_P - [\mu^u_{XUU} P^2 + (2\mu^u_{XXU} - \mu^x_{XXX}) P + (\mu^u_{XU} - 2\mu^x_{XX}) Q \\
&\quad + \mu^x_{XXX}] R_Q,
\end{aligned} \tag{3.20}$$

An admissible cross–section for (3.20) structure equations is given by

$$R_P = R_U = R_X = 0, \quad \text{and} \quad R_Q = 1. \tag{3.21}$$

Substituting (3.21) into (3.20) leads to the following normalizations

$$\mu^x_{XXX} = \mu^u_{XXXX} = \mu^u_{XXXU} = 0, \quad (\text{mod } \omega^x, \omega^u, \omega^p, \omega^q). \tag{3.22}$$

Using the prolongation formula (3.13) for the coefficients of a vector field, we can choose the normalization equations

$$R_{X^i P} = R_{X^i U^j} = 0, \quad i, j \geq 0, \tag{3.23}$$

which cause to normalize the Maurer–Cartan forms $\mu^x_{X^{i+2}}, \mu^u_{X^{i+2} U^j}$.

The group differential of remaining unnormalized second order differential invariants are

$$\begin{aligned}
d_\mathcal{G} R_{QQ} &\equiv -(3\mu^x_X + \mu^u_U) R_{QQ}, \\
d_\mathcal{G} R_{PQ} &\equiv 3\mu^u_{UU} + \mu^u_{XU} - \mu^u_U R_{PQ} - \mu^u_{XU} R_{QQ}, \\
d_\mathcal{G} R_{UQ} &\equiv \mu^u_{UU} - 3\mu^x_X R_{UQ} - \mu^u_{XU} R_{PQ}, \\
d_\mathcal{G} R_{XQ} &\equiv (\mu^u_U - 4\mu^x_X) R_{XQ}, \\
d_\mathcal{G} R_{UP} &\equiv 3\mu^u_{XXUU} - 2\mu^u_{XUU} + (\mu^u_U - 2\mu^x_X) R_{UP} - \mu^u_{XU} R_{PP} - \mu^u_{XU} R_{UQ}, \\
&\vdots
\end{aligned} \tag{3.24}$$

At this stage, we can put

$$R_{UQ} = R_{UP} = R_{PP} = R_{XP} = 0, \quad \text{and} \quad R_{PQ} = R_{XQ} = 1, \tag{3.25}$$

in (3.20) and normalize the Maurer–Cartan forms

$$\mu^u_{XXUU} = \mu^u_{XUU} = 0, \qquad \mu^x_{XX} = \mu^u_{XU} = \mu^u_{UU}, \qquad (\text{mod } \omega^x, \omega^u, \omega^p, \omega^q). \tag{3.26}$$

Considering the structure equations (3.24) and third order differential invariant

$$d_{\mathcal{G}} R_{XXQ} \equiv \mu^u_{XU} + \mu^u_{XU} R_{XQ} - 3\mu^x_X R_{XXQ} - 2R_{XPQ}\omega^q, \tag{3.27}$$

and using the cross–section (3.25), we obtain

$$\mu^u_{UU} \equiv \frac{1}{3}[\mu^u_U + 2\mu^u_{XU} R_{QQ}], \qquad \mu^u_{XU} \equiv -\frac{1}{2}[3\mu^x_X R_{XXQ} + 2R_{XPQ}\omega^q]. \tag{3.28}$$

To obtain the structure equations of the prolonged coframe $\{\omega^x, \omega^u, \omega^p, \omega^q, \mu^x_X, \mu^u_U, \mu^u_{XU}\}$, we substitute the expressions (3.28) into the structure equations (2.21) and the we have the following results:

$$\begin{aligned}
d\omega^x &\equiv \mu^x_X \wedge \omega^x, \\
d\omega^u &\equiv \omega^x \wedge \omega^p + \mu^u_U \wedge \omega^u, \\
d\omega^p &\equiv \omega^x \wedge \omega^q + \mu^u_{XU} \wedge \omega^u + (\mu^u_U - \mu^x_X) \wedge \omega^p, \\
d\omega^q &\equiv \mu^u_{XU} \wedge \omega^p + \mu^u_U \wedge \omega^q, \\
d\mu^x_X &\equiv \frac{3}{2} R_{XXQ} \mu^x_X \wedge \omega^x - R_{XPQ} \omega^x \wedge \omega^q, \\
d\mu^u_U &\equiv \frac{3}{2} R_{XXQ} \mu^x_X \wedge \omega^x - R_{XPQ} \omega^x \wedge \omega^q - \frac{1}{3} \mu^u_U \wedge \mu^u - \frac{2}{3} R_{QQ} \mu_{XU} \wedge \omega^u, \\
d\mu^u_{XU} &\equiv -\frac{1}{3} \mu^u_U \wedge \omega^p - \frac{2}{3} R_{QQ} \mu^u_{XU} \wedge \omega^p + R_{XPQ} \mu^x_X \wedge \omega^q.
\end{aligned} \tag{3.29}$$

Further normalization of the Maurer–Cartan forms $\mu^x_X$, $\mu^u_U$, $\mu^u_{XU}$ depends on the values of invariants $R_{QQ}$, $R_{XPQ}$, the Wunschmann invariant $R_{XXQ}$. If the three invariants are equal to zero then all higher order invariants are automatically zero. According to [5, 18], the necessary and sufficient conditions for a third–order ODE to be linearizable by a fiber–preserving transformation is that $R_{QQ} = R_{XPQ} = R_{XXQ} = 0$. If certain of the invariants $R_{QQ}$, $R_{XPQ}$, $R_{XXQ}$ are non-zero then some of the Maurer–Cartan forms $\mu^x_X$, $\mu^u_U$, $\mu^u_{XU}$ can be normalized, leading to the different branches of the equivalence problem. In generic case $R_{XPQ} R_{XXQ} \neq 0$, for example, we can normalize all three Maurer–Cartan forms by setting $R_{XPQ} = R_{XXQ} = 1$ and $R_{QQ} = 0$. This yields

$$\begin{aligned}
\mu^x_X &\equiv \frac{1}{6}[R_{XXXQ}\omega^x + R_{XXUQ}\omega^u + R_{XXPQ}\omega^p + R_{XXQQ}\omega^q], \\
\mu^u_U &\equiv -\frac{3}{4}[R_{XXXQ}\omega^x + R_{XXUQ}\omega^u + R_{XXPQ}\omega^p + (R_{XXQQ} + 2)\omega^q], \\
\mu^u_{XU} &\equiv -\frac{1}{4}[R_{XXXQ}\omega^x + R_{XXUQ}\omega^u + R_{XXPQ}\omega^p + (R_{XXQQ} + 2)\omega^q].
\end{aligned} \tag{3.30}$$

Substituting the latter expressions into the structure equations of $\omega^x$, $\omega^u$, $\omega^p$, $\omega^q$ in (3.29) we

The equivariant moving frame method of third order differential equations       11

obtain

$$\begin{aligned}
d\omega^x &\equiv \frac{1}{6}[R_{XXUQ}\omega^u \wedge \omega^x + R_{XXPQ}\omega^p \wedge \omega^x + R_{XXQQ}\omega^q \wedge \omega^x],\\
d\omega^u &\equiv \omega^x \wedge \omega^p - \frac{3}{4}[R_{XXXQ}\omega^x \wedge \omega^u + R_{XXPQ}\omega^p \wedge \omega^u + (R_{XXQQ}+2)\omega^q \wedge \omega^u],\\
d\omega^p &\equiv \omega^x \wedge \omega^q + \frac{7R_{XXUQ}-3R_{XXPQ}}{12}\,\omega^p \wedge \omega^u - \frac{R_{XXXQ}}{4}\,\omega^x \wedge \omega^u \\
&\quad -\frac{R_{XXQQ}+2}{4}\,\omega^q \wedge \omega^u - \frac{7R_{XXXQ}}{12}\,\omega^x \wedge \omega^p - \frac{7R_{XXQQ}+18}{12}\,\omega^q \wedge \omega^p,\\
d\omega^q &\equiv \frac{3R_{XXPQ}-R_{XXQQ}-2}{4}\,\omega^q \wedge \omega^p - \frac{1}{4}[R_{XXXQ}\omega^x \wedge \omega^p + R_{XXUQ}\omega^u \wedge \omega^p]\\
&\quad -\frac{3}{4}[R_{XXXQ}\omega^x \wedge \omega^q + R_{XXUQ}\omega^u \wedge \omega^q].
\end{aligned} \qquad (3.31)$$

## Acknowledgements

It is a pleasure to thank the anonymous referees for their constructive suggestions and helpful comments which have materially improved the presentation of the paper. The authors wish to express their sincere gratitude to Francis Valiquette for his critical comments, useful advise and suggestions.